\documentclass[12pt]{amsart}

\theoremstyle{plain}    
\newtheorem{thm}{Theorem}[section]
\numberwithin{figure}{section} 
\theoremstyle{plain}    
\newtheorem{cor}[thm]{Corollary} 
\newtheorem{lemma}[thm]{Lemma} 
\newtheorem{prop}[thm]{Proposition}
\theoremstyle{remark}    
\newtheorem{claim}{Claim}[thm]
\newtheorem{defi}[thm]{Definition}

\usepackage{amscd,amssymb,doublespace,euscript}
\newcount\theTime
\newcount\theHour
\newcount\theMinute
\newcount\theMinuteTens
\newcount\theScratch
\theTime=\number\time
\theHour=\theTime
\divide\theHour by 60
\theScratch=\theHour
\multiply\theScratch by 60
\theMinute=\theTime
\advance\theMinute by -\theScratch
\theMinuteTens=\theMinute
\divide\theMinuteTens by 10
\theScratch=\theMinuteTens
\multiply\theScratch by 10
\advance\theMinute by -\theScratch

\def\today{{\number\day\space
 \ifcase\month\or
  January\or February\or March\or April\or May\or June\or
  July\or August\or September\or October\or November\or December\fi
 \space\number\year}}

\newcommand\Abar{{\overline A}}

\newcommand\Ad{{\rm Ad}}

\newcommand\Afr{{\mathfrak A}}

\newcommand\Ahat{{\widehat A}}

\newcommand\alg{{\textrm{alg}}}

\newcommand\alphabar{{\overline{\alpha}}}

\newcommand\At{{\widetilde A}}

\newcommand\at{{\tilde a}}

\newcommand\Aut{{\rm Aut}}

\newcommand\Bbar{{\overline B}}

\newcommand\betabar{{\overline{\beta}}}

\newcommand\betahat{{\hat\beta}}

\newcommand\clspan{{\overline{\mathrm{span}}\,}}

\newcommand\Cpx{{\mathbf C}}

\newcommand\dif{\mbox{\it d}}

\newcommand\Ec{{\mathcal{E}}}

\newcommand\eqdef{{\;\overset{\mbox{\scriptsize def}}{=}}}

\newcommand\Fc{{\mathcal{F}}}

\newcommand\freeprod{\operatornamewithlimits{\ast}}

\newcommand\freeprodi{{\operatornamewithlimits{\ast}_{\iota\in I}}}

\newcommand\HEu{{\EuScript H}}                   

\newcommand\Ht{{\widetilde H}}

\newcommand\Ints{{\mathbf Z}}

\newcommand\KEu{{\EuScript K}}                   

\newcommand\LEu{{\EuScript L}}                   

\newcommand\lspan{\mathrm{span}\,}

\newcommand\nm[1]{\|#1\|}

\newcommand\Nats{{\mathbf N}}

\newcommand\Oc{{\mathcal{O}}}

\newcommand\phit{{\tilde\phi}}

\newcommand\pihat{{\hat\pi}}

\newcommand\psit{{\tilde\psi}}

\newcommand\rank{\mathrm{rank}\,}

\newcommand\Reals{{\mathbf R}}

\newcommand\Tc{{\mathcal{T}}}

\newcommand\Tcirc{{\mathbf T}}

\newcommand\Ut{{\widetilde U}}

\newcommand\VEu{{\EuScript V}}                   

\begin{document}

\renewcommand{\thefootnote}{\fnsymbol{footnote}}

\pagestyle{myheadings}

\title{Exactness of Cuntz--Pimsner C$^*$--algebras}
 
\author{Kenneth J.\ Dykema, \quad Dimitri Shlyakhtenko$^\dag$}

\begin{abstract}\begin{spacing}{1.0}
Let $H$ be a full Hilbert bimodule over a $C^*$-algebra $A$.  We show that the
Cuntz--Pimsner algebra associated to $H$ is exact if and only if $A$ is exact.
Using this result, we give alternative proofs for exactness of reduced
amalgemated free products of exact $C^*$--algebras.  In the case that $A$ is a
finite--dimensional $C^*$--algebra, we also show that the Brown--Voiculescu
topological entropy of Bogljubov automorphisms of the Cuntz--Pimsner algebra
associated to an $A,A$ Hilbert bimodule is zero.
\end{spacing}\end{abstract}

\maketitle

\markboth{}{}
 
\hfil 1 October, 1999 \hfil

\begin{spacing}{1.2}

\footnotetext[2]{{\rm Supported by an NSF postdoctoral fellowship.}}

\section*{Introduction and description of results.}

A C$^*$--algebra $A$ is said to be {\em exact} if the functor $B\mapsto B\otimes_{\min}A$
preserves short exact sequences of C$^*$--algebras and $*$--homomorphisms.
Recently there has been great progress in understanding exact C$^*$--algebras, much of it due to
E.~Kirchberg and to Kirchberg and S.\ Wassermann; (see~\cite{Ki:nonss}, \cite{Kirchberg:ComUHF},
\cite{Ki95}, \cite{KW:exgp}, \cite{KW:perm}).
A good general reference for exact C$^*$--algebras is Wassermann's monograph~\cite{Wa:SNU}.

In~\cite{Pimsner:C-P}, M.~Pimnser introduced a construction of C$^*$--algebras $E(H)$ and $O(H)$,
respectively called the {\it extended Cuntz--Pimsner algebra} and the {\it Cuntz--Pimsner algebra} of a Hilbert
C$^*$--bimodule $H$ over a C$^*$--algebra $B$; we will review his
construction at the beginning of~\S\ref{sec:ExCP}.
(In this paper we always assume that the Hilbert module $H$ is full as a right $B$--module, i.e.\ that
$\clspan\{\langle \xi,\eta\rangle\mid\xi,\eta\in H\}=B$;
for a good general refernce on Hilbert C$^*$--modules, see the monograph of Lance~\cite{Lance}.)
The C$^*$--algebras $E(H)$ and $O(H)$ are quite important and appear in several areas of operator theory.
They are central in work of Muhly and Solel~\cite{MS:tens}
on triangular operator algebras analogous to the algebra of analytic Toeplitz operators on the circle,
and their ideal structures have been studied in~\cite{DPZ}, \cite{KPW} and~\cite{MS:simpl}, (see also~\cite{Pin:iCP}).

Moreover, the C$^*$--algebra $E(H)$ is related to freeness in the sense of Voiculescu~\cite{V85},
(see also the book~\cite{VDN92}).
For example, Speicher~\cite{Speicher:Mem} has proved that if $H=H_1\oplus H_2$ then $E(H)$ is isomorphic to the
reduced amalgmated free product of C$^*$--algebras $E(H_1)$ and $E(H_2)$, amalgamating over $B$
with respect to the canonical conditional expectations $E(H_\iota)\to B$.
The algebras $E(H)$ are also the natural setting for operator--valued analogues of the Gaussian functor.
The extended Cuntz--Pimsner algebras have been important in work
related to freeness~\cite{shlyakht:amalg}, \cite{shlyakht:semicirc} of the second named author.

In~\S\ref{sec:ExCP} of this paper we show
that if $B$ is an exact C$^*$--algebra and if $H$ is any Hilbert $B,B$--bimodule then
the C$^*$--algebras $E(H)$ and $O(H)$ are exact.
The inspiration for our investigation came from the first named author's recent result~\cite{D:ExactFP} that
every C$^*$--algebra arising as the reduced amalgamated free product of exact C$^*$--algebras is exact.
Moreover, our proof here of exactness of $E(H)$ resembles at its core 
the proof found in~\cite{D:ExactFP} of exactness of reduced amalgamated free products;
in both cases, exactness is proved by in some sense compressing to words of length $n$ and proving exactness
of the resulting compressed C$^*$--algebra by using a chain of ideals of length $n+1$.

This paper's main result on exactness of $E(H)$ can be used to give a new proof of the main result of~\cite{D:ExactFP}.
More sepcifically, in~\S\ref{sec:ExactAFP} we show that if $(A,\phi)=(A_1,\phi_1)*_B(A_2,\phi_2)$ is a reduced amalgamated
free product of C$^*$--algebras (amalgamating over $B$ with respect to conditional expectations $\phi_\iota:A_\iota\to B$)
then $A$ is a quotient of a subalgebra of a quotient of a subalgebra
of $E(H)$, for some Hilbert C$^*$--bimodule over $A_1\oplus A_2$;
if $A_1$ and $A_2$ are exact then from the exactness of $E(H)$ we may conclude that $A$ is exact.
Before giving this argument, however, we consider in~\S\ref{sec:ExactFP} a special case and give an
easier argument showing that if $(A,\phi)=(A_1,\phi_1)*(A_2,\phi_2)$
is a reduced free product (i.e. amalgamating over only the scalars, $\Cpx$), then $A$ can be embedded in $E(H)$ for some Hilbert
C$^*$--bimodule over $A_1\otimes_{\min}A_2$;
this in turn implies the special case of the main result of~\cite{D:ExactFP} that the class of unital exact C$^*$--algebras is
closed under taking reduced free products.
(See~\cite{D:ExGp} for a simpler version of the argument of~\cite{D:ExactFP} in a special case.)

Topolgical entropy for automorphisms of unital nuclear C$^*$--algebras was invented by Voiculescu~\cite{V95}
and was extended by N.P.\ Brown~\cite{Brown:extopent} to apply to automorphisms of exact C$^*$--algebras.
It has many natural properties, and when applied to automorphisms of commutative C$^*$--algebras
gives the usual topological entropy of a homeomorphism.
In~\S\ref{sec:Bog} we examine the topological entropy
of Bogljubov automorphisms of the algebras $E(H)$;
we show (Theorem~\ref{thm:Bog}) that if $H$ is a Hilbert bimodule over $B$ where $B$ is
finite dimensional then every Bogljubov automorphism of $E(H)$ has topological entropy zero.
The fact that we are only able to give results of this sort when $B$ is finite dimensional
parallels the current state of knowledge about the topological entropy of automorphisms of reduced amalgamated
free product C$^*$--algebras; see the results and questions in~\cite{D:TopEnt}.

We have outlined above the contents of the entire paper except for~\S\ref{sec:TopEntCross};
there we collect some results about exactness of crossed product C$^*$--algebras and the topological entropy
of some automorphisms of them.

\section*{Standard Notation.}
\label{sec:stdnot}

We will use the convention $\Nats=\{0,1,2,\ldots\}$.
The words endomorphism, homomorphism, automorphism and representation when
applied to C$^*$--algebras will mean $*$--endomorphism, $*$--homomorphism,
$*$--automorphism and $*$--representation.

\section{Preliminaries on crossed products.}
\label{sec:TopEntCross}

In this section, we will describe and prove some results about topological entropy
in crossed product C$^*$--algebras.
These results are applications of Proposition~2.6 of
Brown and Choda's paper~\cite{BrownChoda:approxent},
which is in turn based on work of Sinclair and
Smith~\cite{SinclairSmith:CBAPx}.
It has come to our attention that
M.\ Choda~\cite{choda:comm} has independently proved Proposition~\ref{prop:htbetahat},
but for completeness we will provide a proof here.

Let us begin by recalling the construction of reduced crossed products,
thereby introducing the notation we will use.
Let $A$ be a C$^*$--algebra with a faithful and nondegenerate representation $\sigma:A\to B(\HEu)$,
where $\HEu$ is a Hilbert space.
Let $G$ be a group taken with discrete topology and let $G\ni g\mapsto\alpha_g\in\Aut(A)$ be
an action of $G$ on $A$ via automorphisms.
Let $\pi:A\to B(\ell^2(G,\HEu))$ be the representation given by
$$ \bigl(\pi(a)\xi)(h)=\sigma\bigl(\alpha_{h^{-1}}(a)\bigr)\bigl(\xi(h)\bigr),
\qquad(a\in A,\,\xi\in\ell^2(G,\HEu),\,h\in G) $$
and let $\lambda$ be the unitary representation of $G$ on $\ell^2(G,\HEu)$ given by
$$ (\lambda_g\xi)(h)=\xi(g^{-1}h),\qquad(\xi\in\ell^2(G,\HEu),\,g,h\in G). $$
We then have
$$ \lambda_{g^{-1}}\pi(a)\lambda_g=\pi\bigl(\alpha_g(a)\bigr),\qquad(a\in A,\,g\in G) $$
and the {\em reduced crossed product C$^*$--algebra} $\Ahat=A\rtimes_\alpha G$ is the norm closure of
the linear span of $\{\pi(a)\lambda_g\mid a\in A,\,g\in G\}$.
It is known that the C$^*$--algebra $\Ahat$ is independent of the choice of $\sigma$.

\begin{lemma}
\label{lemma:betahat}
Let $\Ahat=A\rtimes_\alpha G$ be a reduced crossed product C$^*$--algebra, as described above.
Suppose $\beta\in\Aut(A)$ and $\beta$ commutes with $\alpha_g$ for every $g\in G$.
Then there is a unique automorphism $\betahat\in\Aut(\Ahat)$ such that
$\betahat(\pi(a)\lambda_g)=\pi(\beta(a))\lambda_g$ for every $a\in A$ and $g\in G$.
\end{lemma}

\begin{proof}
Uniqueness is clear.
In the notation used above, we may without loss of generality take the representation $\sigma:A\to B(\HEu)$
to be so that the automorphism $\beta$ is spatially implemented, i.e.\ so that there is a unitary $V\in B(\HEu)$
with $V^*\sigma(a)V=\sigma(\beta(a))$ for every $a\in A$.
Then equating $\ell^2(G,\HEu)$ with $\ell^2(G)\otimes\HEu$ we have the unitary $1\otimes V$ and
$(1\otimes V)^*\pi(a)\lambda_g(1\otimes V)=\pi(\beta(a))\lambda_g$ for every $a\in A$ and $g\in G$.
Let $\betahat=\Ad_{1\otimes V}$.
\end{proof}

For the following theorem,
let us note that the C$^*$--algebra crossed product $A\rtimes_\alpha G$
of an exact C$^*$--algebra $A$ by an action of an amenable countable group $G$ is exact
by~\cite[Proposition 7.1]{Kirchberg:ComUHF}.

\begin{prop}
\label{prop:htbetahat}
Let $A$ be an exact C$^*$--algebra, let $G$ be an amenable countable group,
let $g\mapsto\alpha_g$ be an acton of $G$ on $A$ via automorphisms and let
$\Ahat=A\rtimes_\alpha G$ be the (reduced) crossed product C$^*$--algebra.
Suppose that $\beta\in\Aut(A)$ and that $\beta$ commutes with $\alpha_g$ for every $g\in G$.
Let $\betahat\in\Aut(\Ahat)$ be the automorphism found in Lemma~\ref{lemma:betahat}.
Then $ht(\betahat)=ht(\beta)$.
\end{prop}

\begin{proof}
Because $\betahat\circ\pi=\pi\circ\beta$, using Brown's result~\cite[Proposition 2.1]{Brown:extopent}
that $ht$ is monotone, we have $ht(\betahat)\ge ht(\beta)$.
We will use the notation from the beginning of this section and we will
denote by $\pihat:\Ahat\to B(\ell^2(G,\HEu))$ the inclusion arising from the construction.
In order to show that $ht(\betahat)\le ht(\beta)$, it will suffice to show that
$ht(\pihat,\betahat,\omega,\delta)\le ht(\beta)$ for every $\delta>0$ and for every finite subset $\omega$ of
$\{\pi(a)\lambda_g\mid a\in A,\,g\in G\}$.
Let $\nu$ and $K$ be finite subsets of $A$ and respectively $G$, so that
$\omega\subseteq\{\pi(a)\lambda_g\mid a\in\nu,\,g\in K\}$.
Let $\eta>0$ and let $F$ be a finite subset of $G$ so that $|F\cap gF|\ge(1-\eta)|F|$ for every $g\in K$;
($F$ exists by amenability of $G$).
Let $\nu'=\{\alpha_{t^{-1}}(a)\mid a\in\nu,\,t\in F\}$, let $n$ be a positive integer and let
$k=rcp(\sigma,\nu'\cup\beta(\nu')\cup\cdots\cup\beta^{n-1}(\nu'),\eta)$.
Let $\phi:A\to M_k(\Cpx)$ and $\psi:M_k(\Cpx)\to B(\HEu)$ be completely positive contractions
such that for every $a\in\nu$  and every $j\in\{0,1,\ldots,n-1\}$,
$\nm{\psi\circ\phi\bigl(\beta^j(a)\bigr)-\sigma\bigl(\beta^j(a)\bigr)}<\eta$.
Let $f=|F|^{-1/2}1_F\in\ell^\infty(G)$ be the normalized characteristic function of $F$.
Let $\Phi:\Ahat\to M_{|F|}(\Cpx)\otimes M_k(\Cpx)$ and $\Psi:M_{|F|}(\Cpx)\otimes M_k(\Cpx)\to B(\ell^2(G,\HEu))$
be the completely positive contractions defined in~\cite[Proposition 2.5]{BrownChoda:approxent}.
By~\cite[Proposition 2.6]{BrownChoda:approxent} we have, for every $\at\in A$ and $g\in G$,
$$ \Psi\circ\Phi\bigl(\pi(\at)\lambda_g\bigr)
=\frac1{|F|}\sum_{t\in F\cap gF}\pi\bigl(\alpha_t\circ\psi\circ\phi\circ\alpha_{t^{-1}}(\at)\bigr)\lambda_g. $$
If $g\in K$ and $\at=\beta^j(a)$ for some $0\le j\le n-1$ and $a\in\nu$ then for every
$t\in F\cap gF$ we have
$$ \nm{\alpha_t\circ\psi\circ\phi\circ\alpha_{t^{-1}}(\at)-\at}
=\nm{\psi\circ\phi\circ\beta^j\circ\alpha_{t^{-1}}(a)-\beta^j\circ\alpha_{t^{-1}}(a)}<\eta $$
because $\alpha_t^{-1}(a)\in\nu'$.
Using that $(1-\eta)|F|\le|F\cap gF|\le|F|$ we obtain
$$ \nm{\Psi\circ\Phi\bigl(\pi(\at)\lambda_g\bigr)-\pi(\at)\lambda_g}
<\eta\bigl(\nm a+1\bigr). $$
We could have chosen $\eta$ so small that $\eta\bigl(\nm a+1\bigr)<\delta$ for every $a\in\nu$,
which would have given the estimate
$$ rcp(\pihat,\omega\cup\betahat(\omega)\cup\cdots\cup\betahat^{n-1}(\omega),\delta)
\le|F|\,rcp(\sigma,\nu'\cup\beta(\nu')\cup\cdots\cup\beta^{n-1}(\nu'),\eta). $$
Therefore, $ht(\pihat,\betahat,\omega,\delta)\le ht(\sigma,\beta,\nu',\eta)\le ht(\beta)$.
\end{proof}

We now turn to the crossed product $A\rtimes_\alpha\Nats$
of a C$^*$--algebra $A$ by a single endomorphism $\alpha$;
this construction was introduced by Cuntz~\cite{Cuntz:On}, when he described his algebras $\Oc_n$
as crossed products of UHF algebras by endormorphisms.
Later~\cite[p.\ 101]{Cuntz:Internal} he pointed out that this construction applies more generally.
See Stacey~\cite{Stacey:xendom} for a more detailed discussion, including the nonunital case,
(we consider only his multiplicity one crossed product).
If $A$ is a C$^*$--algebra and if $\alpha$ is an injective endomorphism of $A$,
let $\Abar$ be the inductive limit of the system $A\overset\alpha\to A\overset\alpha\to\cdots$,
with corresponding injective homomorphisms $\mu_n:A\to\Abar$, ($n\in\Nats$).
Let $p$ denote the element $\mu_0(1)$ of $\Abar$ if $A$ is unital, and the corresponding element of
the multiplier algebra of $\Abar$ if $A$ is nonunital.
There is an automorphism $\alphabar$ of $\Abar$ given by $\alphabar(\mu_n(a))=\mu_n(\alpha(a))$,
with inverse $\mu_n(a)\mapsto\mu_{n+1}(a)$.
Then the crossed product $\Ahat=A\rtimes_\alpha\Nats$ is defined to be the hereditary C$^*$--subalgebra
$p\bigl(\Abar\rtimes_\alphabar\Ints\bigr)p$ of the crossed product of $\Abar$ by $\alphabar$.
The map $\mu_0$ followed by the embedding of $\Abar$ into $\Abar\rtimes_\alphabar\Ints$ gives
an embedding $\pi:A\to\Ahat$, and the compression by $p$ of the unitary in
$\Abar\rtimes_\alphabar\Ints$ implementing $\alphabar$ is an isometry $S$ belonging to $\Ahat$
if $A$ is unital and to the multiplier algebra of $\Ahat$ if $A$ is nonunital,
and satisfying
\begin{equation}
\label{eq:SaS}
S\pi(a)S^*=\pi(\alpha(a)),\quad(a\in A).
\end{equation}
If $A$ is unital then $\Ahat$ is the universal unital C$^*$--algebra
generated by a copy $\pi(A)$ of $A$ and an isometry $S$ satisfying~(\ref{eq:SaS});
if $A$ is nonunital then $\Ahat$ satisfies a similar universal property and
is the closed linear span of the set of all elements of the forms
$\pi(a)S^k$ and $(S^*)^k\pi(a)$ for $k\ge0$ and $a\in A$;
see~\cite{Stacey:xendom}.

\begin{lemma}
\label{lem:crossprod}
Let $A$ be an exact C$^*$--algebra and let $\alpha$ be an injective endomorphism of $A$.
Then the crossed product C$^*$--algebra $A\rtimes_\alpha\Nats$ is exact.
\end{lemma}

\begin{proof}
The C$^*$--algebra $\Abar$, being an inductive limit of exact C$^*$--algebras, is exact.
Now~\cite[Proposition 7.1]{Kirchberg:ComUHF} implies that $\Abar\rtimes_\alphabar\Ints$ is exact,
hence that $A\rtimes_\sigma\Nats$ is exact.
\end{proof}

The following lemma follows easily from the universal property for the crossed product by an endormorphism,
but we will exhibit the automrphism $\betahat$ directly, for use in the next proposition.
\begin{lemma}
\label{lemma:betahatend}
Let $A$ be a C$^*$--algebra with an injective endormorphism $\alpha$
and an automorphism $\beta$ that commutes with $\alpha$.
Then there is a unique autormorphism $\betahat$ of $A\rtimes_\alpha\Nats$ satisfying
$\betahat(\pi(a)S^k)=\pi(\beta(a))S^k$ for every $a\in A$ and $k\ge0$.
\end{lemma}

\begin{proof}
Uniqueness is clear.
Let $\Ahat=A\rtimes_\alpha\Nats$;
then $\Ahat=p(\Abar\rtimes_\alphabar\Ints)p$ as above.
The commuting diagram
\[
\begin{CD}
A@>\alpha>>A@>\alpha>>\cdots \\
@VV\beta V @VV\beta V \\
A@>\alpha>>A@>\alpha>>\cdots
\end{CD}
\]
gives rise to an automorphism $\betabar$ of $\Abar$ that commutes with $\alphabar$.
Let $\gamma$ be the automorphism of $\Abar\rtimes_\alphabar\Ints$ arising from $\betabar$
via Lemma~\ref{lemma:betahat}.
Then $\gamma(\Ahat)=\Ahat$ and the restriction of $\gamma$ to $\Ahat$
is the desired automorphism $\betahat$.
\end{proof}

\begin{prop}
\label{prop:htbetahatend}
Let $A$ be an exact C$^*$--algebra with an injective endomorphism $\alpha$
and an automorphism $\beta$ that commutes with $\alpha$;
let $\Ahat=A\rtimes_\alpha\Nats$.
Let $\betahat\in\Aut(\Ahat)$ be the automorphism found in Lemma~\ref{lemma:betahatend}.
Then $ht(\betahat)=ht(\beta)$.
\end{prop}

\begin{proof}
Let us use the notation of the proof of Lemma~\ref{lemma:betahatend}.
Then we have
\begin{equation*}
ht(\beta)\le ht(\betahat)\le ht(\gamma)=ht(\betabar),
\end{equation*}
where the inequalities follow from monotonicity of $ht$
and the equality follows from Proposition~\ref{prop:htbetahat}.
However, Brown's result~\cite[Proposition 2.14]{Brown:extopent} on inductive limit automorphisms
gives $ht(\betabar)=ht(\beta)$.
\end{proof}

\section{Exactness of the Cuntz-Pimsner algebras.}
\label{sec:ExCP}

In this section we prove that the extended Cuntz--Pimsner algebra $E(H)$
and the Cuntz--Pimsner algebra $O(H)$ of a Hilbert $B,B$--bimodule $H$
are exact C$^*$--algebras whenever $B$ is an exact C$^*$--algebra.
We begin by reviewing Pimsner's construction~\cite{Pimsner:C-P} of these algebras
and some facts about them.

Let \( B \) be a C$^*$--algebra and let \( H \) be a Hilbert bimodule
over \( B \).
By this we mean that $H$ is a right Hilbert $B$--module with an injective homomorphism $B\to\LEu(H)$;
we further assume that $\{\langle h_1,h_2\rangle_B\mid h_1,h_2\in H\}$ generates $B$ as a C$^*$--algebra,
where $\langle\cdot,\cdot\rangle_B$ is the $B$--valued inner product on $H$.
Let \( \mathcal{F}(H)=B\oplus \bigoplus _{n\geq 1}H^{(\otimes _{B})n} \)
be the full Fock space over \( H \);
here $H^{(\otimes_B)n}$ denotes the $n$--fold tensor product $H\otimes_BH\otimes_B\cdots\otimes_BH$.
Note that $\Fc(H)$ is a Hilbert $B,B$--bimodule.
For each vector \( h\in H \), the operator
\( l(h):\mathcal{F}(H)\to \mathcal{F}(H) \) defined by
\begin{eqnarray*}
l(h)h_{1}\otimes \dots \otimes h_{n} & = & h\otimes h_{1}\otimes \dots \otimes h_{n},\quad h,h_{1},\dots ,h_{n}\in H\\
l(h)b & = & hb,\quad h\in H,\, b\in B
\end{eqnarray*}
is a bounded adjointable operator on \( \Fc(H) \).
These $l(h)$ are called creation operators and satisfy the relations
\begin{eqnarray}
l(h)^{*}l(g) & = & \langle h,g\rangle _{B},\quad h,g\in H\label{eqn:ls} \\
b_{1}l(h)b_{2} & = & l(b_{1}hb_{2}),\quad h\in H,\, b_{1},b_{2}\in B.\label{eqn:ls1} 
\end{eqnarray}

Pimsner defined the \emph{extended Cuntz-Pimsner algebra} \( E(H)\subset\LEu(\Fc(H)) \) to be
\[
E(H)=C^{*}(l(h):h\in H).\]
(where his notation is $\Tc_H$).
Since we assumed that $B$ is generated by the set of inner products $\langle h_1,h_2\rangle$,
the copy of $B$ acting on the left of $\Fc(H)$ is contained in $E(H)$.
Pimsner showed~\cite[Theorem 3.4]{Pimsner:C-P} that \( E(H) \) is in fact
the universal C$^*$--algebra generated by \( B \) and
elements \( l(h) \), satisfying relations (\ref{eqn:ls}) and (\ref{eqn:ls1}).
The orthogonal projection onto \( B\subset \mathcal{F}(H) \) defines a canonical
conditional expectation \( \mathcal{E} \) from \( E(H) \) onto \( B \).

If \( K\subset H \) is a Hilbert subbimodule, then \( C^{*}(l(h):h\in K)\cong E(K) \),
so that there is an inclusion \( E(K)\subset E(H) \);
this inclusion preserves \( \mathcal{E} \).
Note furthermore that if $H'$ is a closed $\Cpx$--linear subspace of $H$ and
if $B'$ is a C$^*$--subalgebra of $B$ such that
\begin{eqnarray*}
\langle h_1,h_2\rangle_B\in B'&\quad(h_1,h_2\in H') \\
b_1hb_2\in H'&\quad(b_1,b_2\in B',\,h\in H')
\end{eqnarray*}
then $H'$ is a Hilbert bimodule over $B'$ and $E(H')\subset E(H)$.

\begin{thm}
\label{thm:EHexact}
Let $B$ be a C$^*$--algebra and let \( H \) be a Hilbert $B,B$--bimodule
such that $\{\langle h_1,h_2\rangle\mid h_1,h_2\in H\}$ generates $B$.
Then \( E(H) \) is exact if and only if \( B \) is exact.
\end{thm}

\begin{proof}
Since C$^*$--subalgebras of exact C$^*$--algebras are exact and since
\( B\subset E(H) \),  if \( E(H) \) is exact then \( B \) is exact.

Assume now that \( B \) is exact, and let us show that $E(H)$ is exact.
There is a net, ordered by inclusions, of pairs $(B'_\lambda,H'_\lambda)$
where each $B'_\lambda$ is a separable C$^*$--subalgebra of $B$
and where each $H'_\lambda$ is a seperable closed linear subspace of $H$
such that the restiction of the usual operations makes $H'_\lambda$ a Hilbert bimodule over $B'_\lambda$
and such that $\overline{\bigcup_\lambda B'_\lambda}=B$ and $\overline{\bigcup_\lambda H'_\lambda}=H$.
From the inclusions $E(H_\lambda')\subset E(H)$ mentioned early in this section,
we see that $E(H)$ is the direct limit of the $E(H'_\lambda)$.
Hence we may and do assume without loss of generality that $B$ and $H$ are separable.

Let \( \Ht=H\oplus B \).
Since \( E(H)\subset E(\Ht) \), it will be sufficient to prove that \( E(\Ht) \) is exact.
Denote by \( \xi \in \Ht \)
the vector \( 0\oplus 1_{B} \), and let \( L=l(\xi ) \). Then \( L \) satisfies
the following relations: 
\begin{eqnarray*}
L^{*}L & = & 1,\\
L^{*}l(h) & = & 0,\quad h\in H\\
l(h)^{*}L & = & 0,\quad h\in H.
\end{eqnarray*}
The C$^*$--algebra $E(\Ht)$ is the closed linear span of the set of all
elements of the form
\begin{equation}
\label{eqn:Wword}
W=b_{0}l(h_{1})^{g(1)}b_{1}l(h_{2})^{g(2)}b_2\cdots l(h_n)^{g(n)}b_{n}
\end{equation}
where \( n\geq 0 \), $b_j\in B$, \( g(j)\in \{*,\cdot \} \), \( h_{1},\ldots ,h_{n}\in \Ht \)
and where $l(h_j)^{g(j)}=l(h_j)$ if $g(j)=\cdot$.

Consider the unitary \( \alpha _{t}:\Ht\to \Ht \) given by \( \alpha _{t}(h)=e^{2\pi it}h \),
$(h\in \Ht$).
Denote by \( \beta _{t} \) the resulting automorphism \( E(\alpha _{t}) \)
of \( E(\Ht) \).
Thus \( \beta _{t}(L)=e^{2\pi it}L \), \( \beta _{t}(l(h))=e^{2\pi it}l(h) \),
($h\in H$) and $\beta_t(b)=b$, ($b\in B)$.
Note that $t\mapsto\beta_t$ is an action of the group $\Tcirc=\Reals/\Ints$ on $E(H)$.
Let \( A \) be the fixed point subalgebra of $\beta$, i.e.\
$A=\{a\in E(H)\mid \forall t\in\Tcirc\,\beta_t(a)=a\}$.

\begin{claim}
\label{claim:Aspan}
\( A \) is the closed linear span of the set of operators of the form~\eqref{eqn:Wword}
for which \(\#\{i:g(i)=*\}=\#\{i:g(i)=\cdot \} \).
\end{claim}
\begin{proof}
If $W$ is of the form~\eqref{eqn:Wword} then
\( \beta _{t}(W)=e^{2\pi i(\#\{i:g(i)=\cdot \}-\#\{i:g(i)=*\})t}W \).
Hence if \( \#\{i:g(i)=*\}=\#\{i:g(i)=\cdot \} \) then \( W\in A \).
The map
\( \Phi (T)=\int_0^1\beta _{t}(T)dt \) is a faithful conditional
expectation from \( E(\Ht) \) onto \( A \);
letting \( W \) be as above, if \( \#\{i:g(i)=*\}=\#\{i:g(i)=\cdot \} \) then $\Phi(W)=W$
while otherwise $\Phi(W)=0$.
If \( T\in A \), then \( T \) can be approximated by
a linear combination of operators of the form~\eqref{eqn:Wword}.
Since \( \Phi (T)=T \) by assumption, it then follows that \( T \) can be approximated
by a linear span of operators of the form~\eqref{eqn:Wword}
for which \( \#\{i:g(i)=*\}=\#\{i:g(i)=\cdot \} \).
This proves Claim~\ref{claim:Aspan}.
\end{proof}

\begin{claim}
\label{claim:EHgen}
\( E(\Ht)=C^{*}(A,L) \).
\end{claim}
\begin{proof}
It is sufficient to show that any operator $W$ of the form~(\ref{eqn:Wword})
can be written as \( W=(L^{*})^{k}W' \) or \( W=W'L^{k} \) for some \( W'\in E(H) \)
and \( k\in \Nats \).
Let \( k=\#\{i:g(i)=\cdot \}-\#\{i:g(i)=*\} \).
If \( k=0 \), then \( W\in A \) and we are done.
If \( k>0 \) then since \( L^{*}L=1 \) we have
\( W=(W(L^{*})^{k})L^{k} \) and \( W'=W(L^{*})^{k}\in A \).
If \( k<0 \), then \( W=(L^{*})^{k}L^{k}W \), and \( W'=L^{k}W\in A \).
This proves Claim~\ref{claim:EHgen}.
\end{proof}

\begin{claim}
\label{claim:crossed}
\( \Psi :a\mapsto LaL^{*} \) defines an injective endomorphism of \( A \).
\( E(\Ht) \) is isomorphic to the (universal) crossed product of \( A \) by this endomorphism, namely
$E(\Ht)\cong A\rtimes_\Psi\Nats$.
\end{claim}
\begin{proof}
\( \Psi  \) is an injective endomorphism because \( L^{*}L=1 \) and $L\ne0$.
Let $C=A\rtimes_\Psi\Nats$ and let \( V\in C \) denote
the isometry arising from the crossed product contruction and implementing \( \Psi  \);
thus we have \( VaV^{*}=\Psi (a) \), \( a\in A \).
As is well known, there exists a continuous family $(\gamma_t)_{t\in\Tcirc}$ of automorphisms \( \gamma _{t} \) of \( C \),
such that \( \gamma _{t}(V)=e^{2\pi it}V \), and \( \gamma _{t}(a)=a \) for every
\( a\in A \), (where we identify the circle $\Tcirc$ with $\Reals/\Ints$).
Let \( \Gamma :C\to A \) be the conditional expectation \( \Gamma (T)=\int_0^1\gamma _{t}(T)dt \).
Then \( \Gamma  \) is faithful. By universality of \( C \), there is a surjective
map \( \rho :C\to E(\Ht) \), such that \( \rho (a)=a \) for \( a\in A \), and
\( \rho(V)=L \).
Let \( T\in \ker \rho  \).
Then \( T^{*}T\in \ker \rho  \).
Since \( \beta \circ \rho =\rho \circ \gamma  \), we have \( \gamma _{t}(T^{*}T)\in \ker \rho  \),
so that \( \Gamma (T^*T)\in \ker \rho  \).
But \( \rho |_{A} \) is injective,
so \( \Gamma (T^*T)=0 \) and hence \( T^*T=0 \) and $T=0$.
It follows that \( C\cong E(\Ht) \). 
This proves Claim~\ref{claim:crossed}.
\end{proof}

\begin{claim}
\label{claim:Aexact}
The C$^*$--algebra \( A \) is exact.
\end{claim}
\begin{proof}
Denote by \( A_{n} \) the subspace
of \( A \) that is the closed linear span of the set of words of the form
\( W=b_{0}l(h_{1})^{g(1)}b_{1}l(h_{2})^{g(2)}b_2\cdots l(h_{2m})b_{2m} \)
with \( m\leq n \), and for which \( \#\{i:g(i)=\cdot \}=\#\{i:g(i)=*\} \).
Note that, in light of equations~(\ref{eqn:ls}) and~(\ref{eqn:ls1}),
we may without loss of generality assume that in $W$,
$g(1)=g(2)=\cdots=g(m)=\cdot$ and $g(m+1)=g(m+2)=\cdots=g(2m)=*$.
Now it is easily seen that \( A_{n} \) is a C$^*$--subalgebra of \( A \) and that
\( B=A_{0}\subset \cdots \subset A_{n}\subset A_{n+1}\subset A \)
is an increasing sequence of subalgebras with \(  \bigcup_{n\ge1}A_{n} \) dense in \( A \).
Hence it will suffice to prove that each \( A_{n} \) is exact.

Denote by \( \pi _{n} \) the restriction and compression of the representation of \( A_{n} \)
on \( \mathcal{F}(\Ht) \) to \( \mathcal{F}_{n}(\Ht)=B\oplus \bigoplus _{k\leq n}\Ht^{(\otimes _{B})k} \).
Since in the decomposition
\begin{equation}
\label{eqn:Fcdecomp}
\mathcal{F}(\Ht)=
\bigoplus_{k\geq 0} \bigoplus_{j=0}^n\Ht^{\otimes_B(k(n+1)+j)}=
\bigoplus_{j=0}^n\bigoplus_{k\geq 0}\Ht^{\otimes_B j}\otimes 
\left(\Ht^{\otimes_B(n+1)}\right)^{\otimes_B k} =
\mathcal{F}_n(\Ht)\otimes\mathcal{F}(\Ht^{\otimes_B(n+1)})
\end{equation}
\( A_{n} \) acts on $\Fc(\Ht)$  as \( \pi _{n}(A_n)\otimes 1 \), it follows that \( \pi _{n} \)
is a faithful representation.
Let \( I_{n}\subset A_{n} \) be the closed linear span of the set of all words of the form
\( b_{0}l(h_{1})b_{1}l(h_{2})b_2\cdots l(h_{n})b_nl(h_{n+1})^{*}b_{n+1}\cdots l(h_{2n})^{*}b_{2n} \).
Using the relations~(\ref{eqn:ls}) and~(\ref{eqn:ls1}) one easily sees that $I_n$
is a closed two sided ideal of $A_n$.
Moreover, observing the action of $\pi_n(I_n)$ on $\Fc_n(\Ht)$
one easily sees that,
with respect to  the decomposition $\Fc_n(\Ht)=\Fc_{n-1}(\Ht)\oplus \Ht^{(\otimes_b)n}$,
we have $\pi_n(I_n)=0_{\Fc_{n-1}(\Ht)}+\KEu(\Ht^{(\otimes_B)n})$,
where $\KEu(\Ht^{(\otimes_B)n})$ denotes the algebra of compact operators
on the Hilbert \( B,B \)--bimodule $\Ht^{(\otimes_B)n}$.
Now the quotient \( A_{n}/I_{n} \) is canonically identified with the closed linear span of
all words \( b_{0}l(h_{1})^{g(1)}b_{1}l(h_{2})^{g(2)}\cdots l(h_{2m})^{g(2m)}b_{2m} \) in \( A_{n} \)
for which \( m<n \).
Thus \( A_{n}/I_{n} \) is isomorphic to \( A_{n-1} \), and the canonical inclusion $A_{n-1}\hookrightarrow A_n$
provides a splitting for the short exact sequence \( 0\to I_{n}\to A_{n}\to A_{n}/I_{n}\to 0 \).

We now show exactness of $A_n$ by induction on \( n \).
For \( n=0 \), \( A_{0}\cong B \) is exact by assumption.
Having restricted to the separable case, we have that $\Ht^{(\otimes_B)n}$ is separable
and hence $\KEu(\Ht^{(\otimes_B)n})$ is an exact C$^*$--algebra
by the Kasparov stabilization lemma~\cite{Kasp}. 
Using the induction hypothesis that \( A_{n-1} \) is exact, we find that
in the split exact sequence \( 0\to I_{n}\to A_{n}\to A_{n}/I_{n}\to 0 \) the
algebras \( A_{n}/I_{n}\cong A_{n-1} \) and \( I_{n}\cong\KEu(\Ht^{(\otimes_B)n}) \)
are exact.
By \cite[Proposition 2]{HRV:exact}, (see also~\cite[Proposition 7.1]{Kirchberg:ComUHF}),
\( A_{n} \) is exact.
This completes the proof of Claim~\ref{claim:Aexact}.
\end{proof}

Now we can finish the proof of the theorem.
By Claims~\ref{claim:crossed} and~\ref{claim:Aexact}, $E(\Ht)$ is isomorphic to
the universal crossed product $A\rtimes_\Psi\Nats$ of an exact C$^*$--algebra $A$
by an injective endomorphism $\Psi$.
By Lemma~\ref{lem:crossprod}, $A\rtimes_\Psi\Nats$ is exact.
\end{proof}

\begin{cor}
Let $B$ be a C$^*$--algebra and let $H$ be a Hilbert \( B \),\( B \)--bimodule.
Then the Cuntz-Pimsner algebra $\Oc_H$ associated to $H$ is exact if and only if
\( B \) is exact.
\end{cor}
\begin{proof}
The Cuntz-Pimsner algebra \( \mathcal{O}_{H} \) is defined
as a certain quotient of \( E(H) \) (see \cite{Pimsner:C-P} for details).
If $B$ is exact then $E(H)$ is exact, and it was proved by Kirchberg~\cite{Kirchberg:ComUHF}
that quotients of exact C$^*$--algebras are exact.
Since $\Oc_H$ contains a copy of $B$ as a C$^*$--subalgebra, exactness of $\Oc_H$ implies exactness of $B$.
\end{proof}

Given a C$^*$--algebra \( B \) and a completely-positive map \( \eta :B\to B\otimes M_{n\times n}(\Cpx) \)
the C$^*$--algebra \( \hat{\Phi }(B,\eta ) \) was constructed in~\cite{shlyakht:semicirc}.
This algebra is a subalgebra of \( E(H) \),
where \( H \) is the Hilbert \( B,B \) bimodule associated to \( \eta  \)
(see \cite{rieffel:rigged}). Since \( B\subset \hat{\Phi }(B,\eta ) \), from Theorem~\ref{thm:EHexact}
we immediately have the following.

\begin{cor}
Let \( B \) be a C$^*$--algebra and let \( \eta :B\to B\otimes M_{n\times n} \)
a completely-positive map.
Then the C$^*$--algebra \( \hat{\Phi }(B,\eta ) \) is exact if and only if \( B \) is exact.
\end{cor}

\section{Exactness of Reduced Free Product C$^*$--algebras.}
\label{sec:ExactFP}

Theorem~\ref{thm:EHexact} can be used to give a new proof of the recent
result~\cite{D:ExactFP} that the class of exact unital C$^*$--algebras is closed under
taking reduced amalgamated free products;
this alternative proof will be given in section~\ref{sec:ExactAFP} below.
In this section, however, we give an easier argument
that makes use of Theorem~\ref{thm:EHexact} and proves a special case, namely that 
every reduced free product (with amalgamation over the scalars) of exact C$^*$--algebras is exact.
(See~\cite{D:ExGp} for an easier version of the argument of~\cite{D:ExactFP} in a special case.)

In light of Example~1.4 of~\cite{DR:proj},
(see also Question~1 of~\cite{BD:emb} and the answers provided),
one must be careful about embeddings of reduced free products of C$^*$--algebras.
Thus we provide a proof of the following lemma.
\begin{lemma}
\label{lem:fptens}
Let $N$ be an integer $\ge2$ or $\infty$ and for every $1\le k<N+1$ let $A_k$ be
a unital C$^*$--algebra having a state $\phi_k$ whose GNS representation is faithful.
Let $(A,\phi)=\freeprod_{k=1}^N(A_k,\phi_k)$
be the reduced free product of C$^*$--algebras.
Let $B=\bigotimes_{k=1}^NA_k$ be the minimal tensor product of C$^*$--algebras
and let $\rho=\otimes_{k=1}^N\phi_k$ be the tensor product state.
Let $D_1$ be a C$^*$--algebra with a state $\psi_1$ whose GNS representation is faithful
and having a unitary $u\in D_1$ such that $\psi_1(u^n)=0$ for every nonzero integer $n$.
Let $(D,\psi)=(D_1,\psi_1)*(B,\rho)$ be the reduced free product of C$^*$--algebras.
Consider the embeddings $\pi_k:A_k\to D$ given by $\pi_k(a)=u^kau^{-k}$.
Then there is an injective homomorphism $\pi:A\to D$
whose restriction to $A_k$ is $\pi_k$, for every $k$ and such that $\psi\circ\pi=\phi$.
\end{lemma}

\begin{proof}
Let us abuse notation by writing $A_k$ for all of the corresponding unital subalgebras
$A_k\subseteq A$, $A_k\subseteq B$ and $A_k\subseteq D$
arising from the free product and tensor product constructions, and similarly for
$B\subseteq D$.
The unitary $u$ generates a copy of $C(\Tcirc)$, the continuous functions on
the circle, on which $\psi_1$ is given by integration with respect to Haar measure.
Thus $(C(\Tcirc),\int\cdot\dif\lambda)\subseteq(D_1,\psi_1)$,
where $\dif\lambda$ is Haar measure.
By the main result of~\cite{BD:emb}, without loss of generality we may and do assume
that $(D_1,\psi_1)=(C(\Tcirc),\int\cdot\dif\lambda)$.
It is easily checked that in $(D,\psi)$ the family $(u^kBu^{-k})_{k\in\Ints}$
is free;
letting $\Bbar$ be the C$^*$--subalgebra of $D$ generated by
$\bigcup_{k\in\Ints}u^kBu^{-k}$,
conjugation by $u$ acts as the free shift on $\Bbar$.
As $\Bbar\cup\{u\}$ generates $D$ and as $\Bbar u^k\subseteq\ker\psi$ for every
nonzero integer $k$, we see that $D\cong\Bbar\rtimes\Ints$ and the GNS representation
of the restriction of $\psi$ to $\Bbar$ is faithful on $\Bbar$.
Therefore, we may use the uniqueness of the free product construction
to see that
$(\Bbar,\psi{\restriction}_\Bbar)\cong
\freeprod_{k=-\infty}^\infty(u^kBu^{-k},\psi{\restriction}_{u^kBu^{-k}})$.
Regarding $A_k\subseteq B$, we have the embeddings $\pi_k:A_k\to u^kA_ku^{-k}\subset u^kBu^{-k}$
and these satisfy $\psi\circ\pi_k=\phi_k$.
Hence by the main result of~\cite{BD:emb} there is an injective homomorphism
$\pi:A\to\Bbar\subset D$ extending each $\pi_k$ and satisfying $\psi\circ\pi=\phi$.
\end{proof}

\begin{prop}
\label{prop:AemEH}
Let $I$ be a set having at least two elements and for every $\iota\in I$ let
$A_\iota$ be a unital C$^*$--algebra and let $\phi_\iota$ be a state on $A_\iota$
whose GNS representation is faithful.
Let $(A,\phi)=\freeprodi(A_\iota,\phi_\iota)$ be the
reduced free product of C$^*$--algebras, let $B=\bigotimes_{\iota\in I}A_\iota$
be the minimal tensor product of C$^*$--algebras and let $\rho=\otimes_{\iota\in I}\phi_\iota$
be the tensor product state.
Then there is a Hilbert $B,B$--bimodule $H$ and an injective homomorphism $\pi:A\to E(H)$
such that letting $\Ec:E(H)\to B$ be the canonical vacuum expectation, we have $\rho\circ\Ec\circ\pi=\phi$.
\end{prop}
\begin{proof}
Let \( H \) be the Hilbert \( B,B \)--bimodule associated to \( \rho  \).
This means that \( H \) is obtained by separation and completion of the algebraic
tensor product \( B\otimes_\alg B \)
with respect to the norm induced by the $B$--valued inner product
\[
\langle a\otimes b,a'\otimes b'\rangle =b^{*}\rho (a^{*}a')b',\quad a,a',b,b'\in B, \]
or, in other notation, $H=L^2(B,\rho)\otimes_\Cpx B$.
Denote by \( \xi \in H \) the vector \( 1\otimes 1 \).
Let $\Ec:E(H)\to B$ be the canonical vacuum expectation given by compression with the projection $\Fc(H)\to B$.
Consider \( C=C^{*}(l(\xi ))\subset E(H) \).
By~\cite{shlyakht:amalg},
the restriction of the conditional expectation \( \mathcal{E}:E(H)\to B \)
to \( C \) is scalar-valued;
we denote this restriction by \( \psi  \).
In fact, as is easily seen, \( C \) is isomorphic to the algebra
of Toeplitz operators generated by the nonunitary isometry $l(\xi)$
and $\psi$ is the state whose support is $1-l(\xi)l(\xi)^*$.
We have by \cite[Theorem 2.3]{shlyakht:amalg} that $E(H)$ is a reduced free product,
\( (E(H),\rho\circ\mathcal{E})\cong (C,\psi )*(B,\rho ) \),
because \( l(\xi ) \) satisfies \( l(\xi )^{*}bl(\xi )=\rho (b) \) for all
\( b\in B \).
The algebra $C$ contains a unitary \( u \) with the property that \( \psi (u^{k})=0 \)
for all \( k\in \Ints\setminus \{0\} \);
for example, $u$ can be obtained by contiunuous functional calculus from the semicircular element $l(\xi)+l(\xi)^*$.
Then by Lemma~\ref{lem:fptens} there is a injective homomorphism $\pi:A\to E(H)$
satisfying that $\rho\circ\Ec\circ\pi=\phi$.
\end{proof}

\begin{cor}[\cite{D:ExactFP}]
\label{cor:ExactFPC}
Let $I$ be a set having at least two elements and for every $\iota\in I$ let
$A_\iota$ be a unital C$^*$--algebra and let $\phi_\iota$ be a state on $A_\iota$
whose GNS representation is faithful.
Let $(A,\phi)=\freeprodi(A_\iota,\phi_\iota)$ be the
reduced free product of C$^*$--algebras.
Then $A$ is exact if and only if every $A_\iota$ is exact.
\end{cor}
\begin{proof}
Since each $A_\iota$ is canonically embedded as a C$^*$--subalgebra of $A$,
exactness of $A$ implies exactness of every $A_\iota$.

Suppose that every $A_\iota$ is exact and let $B=\bigotimes_{\iota\in I}A_\iota$
be the minimal tensor product of C$^*$--algebras.
Then $B$ is exact, as is easily seen from the definition of exactness
and by taking inductive limits if necessary.
By Proposition~\ref{prop:AemEH}, $A$ is isomorphic to a C$^*$--subalgebra of $E(H)$, for some
Hilbert $B,B$--bimodule $H$.
By Theorem~\ref{thm:EHexact} \( E(H) \) is exact, and it thus follows that $A$ is exact.
\end{proof}

\section{Exactness of Reduced Amalgamated Free Product C$^*$--algebras.}
\label{sec:ExactAFP}

In this section, we give an alternative proof, using Theorem~\ref{thm:EHexact}, of the result~\cite{D:ExactFP}
that the class of exact unital C$^*$--algebras is closed under taking reduced amalgamated free products.

\begin{prop}
\label{prop:AtemEH}
Let $B$ be a unital C$^*$--algebra and let $A_1$ and $A_2$ be unital C$^*$--algebras each containing
a copy of $B$ as a unital C$^*$--subalgebra and having a conditional expectation
$\phi_\iota:A_\iota\to B$ whose GNS representation is faithful.
Let
\begin{equation*}
(\At,\phit)=(A_1,\phi_1)*_B(A_2,\phi_2)
\end{equation*}
be the reduced amalgamated free product of C$^*$--algebras and let $A=A_1\oplus A_2$.
Then there is a Hilbert $A,A$--bimodule $H$ such that $\At$ is isomorphic to a C$^*$--quotient of a C$^*$--subalgebra
of a C$^*$--quotient of a C$^*$--subalgebra of $E(H)$.
\end{prop}
\begin{proof}
Let $D=B\oplus B=\{(b_1,b_2)\in A\mid b_1,b_2\in B\}$ and
consider the conditional expectation $\phi=\phi_1\oplus\phi_2:A\to D$
and the completely positive map $\eta:A\to A$ given by
$\eta\bigl((a_1,a_2)\bigr)=(\phi_2(a_2),\phi_1(a_1))$.
Note that $\eta$ takes values in $D$ and that $\eta\circ\phi=\eta$.
We now consider the $\eta$--creation operator $L$;
this construction was introduced by Speicher~\cite{Speicher:Mem} and Pimsner~\cite{Pimsner:C-P},
and proceeds as follows.
We take the Hilbert $A,A$--bimodule $H$ which is obtained by separation and completion of
the algebraic tensor product $A\otimes_\alg A$ equipped with the natural left and right actions of $A$ and with the
inner product
\begin{equation*}
\langle a_1\otimes a_2,a_1'\otimes a_2'\rangle=a_2^*\eta(a_1^*a_1')a_2'.
\end{equation*}
We let $\xi\in H$ be the element corresponding to $1\otimes 1\in A\otimes A$;
taking the Fock
space $\Fc(H)$ we let $L=l(\xi)\in E(H)\in\LEu(\Fc(H))$ be the corresponding creation operator.
As is well known and is easily seen using equations~\eqref{eqn:ls} and~\eqref{eqn:ls1},
$L^*aL=\eta(a)$ for every $a\in A$.

\begin{claim}
\label{claim:LAfree}
Let $\Ec:E(H)\to A$ be the conditional expectation defined by compression with the orthogonal projection
$\Fc(H)\to A$ and let $\psi=\phi\circ\Ec:E(H)\to D$.
Then $\{L,L^*\}$ and $A$ are free with respect to $\psi$.
\end{claim}
\begin{proof}
This claim is proved by applying~\cite[Theorem~2.3]{shlyakht:amalg}.
In fact, from its proof, we see that statements~(i), (ii) and~(b) of that theorem imply
statement~(a) of that theorem.
Hence in order to apply~\cite[Theorem~2.3]{shlyakht:amalg} to show freeness of $A$ and $\{L,L^*\}$,
we need only show
\renewcommand{\labelenumi}{(\roman{enumi})}
\begin{enumerate}
\item[($\alpha$)] $\psi(a_1La_2\cdots La_kLa_1'L^*a_2'\cdots L^*a_{\ell+1}')=0$
     whenever $k,\ell\ge0$, $k+\ell>0$ and $a_j,a_j'\in A$;
\item[($\beta$)] $L^*aL=\eta(\psi(a))$ for every $a\in A$.
\end{enumerate}
Indeed, ($\alpha$) and ($\beta$) together show that $L$ is distributed with respect to $\psi$ as
an $(\eta{\restriction}_D)$--creation operator, showing part~(c) of~\cite[Theorem~2.3]{shlyakht:amalg} holds,
while ($\alpha$) is part~(i) and ($\beta$) is part~(ii) of~\cite[Theorem~2.3]{shlyakht:amalg}.
But ($\alpha$) and ($\beta$) follow from the facts that $L$ is distributed as an $\eta$--creation operator
and $\eta\circ\psi(a)=\eta\circ\phi(a)=\eta(a)$ for $a\in A$.
This finishes the proof of Claim~\ref{claim:LAfree}.
\end{proof}

Now let $P=1-L^2(L^*)^2$ and $W=P(L+L^*)P$.
\begin{claim}
\label{claim:PW}
$P$ is a projection, $W$ is a partial isometry, $W=W^*$ and $W^2=P$.
\end{claim}
\begin{proof}
The creation operator $L$ is easily seen to be an isometry;
thus $P$ is a projection.
Clearly $W=W^*$.
Straightforward computations reveal that
\begin{equation}
\label{eqn:WL}
W=L+L^*-L(L^*)^2-L^2L^*
\end{equation}
and then that $W^2=P$.
So Claim~\ref{claim:PW} is proved.
\end{proof}

\begin{claim}
\label{claim:Wunitary}
Let $\Afr=C^*(\{W\}\cup A)\subseteq E(H)$ and let $\pi$ be the GNS representation of $\Afr$
arising from the conditional expectation $\psi$.
Then $\pi(1-P)=0$;
hence $\pi(W)$ is unitary
\end{claim}
\begin{proof}
We must show that
\begin{equation}
\label{eqn:psiaw}
\psi(a_1W^{q_1}\cdots a_kW^{q_k}a_{k+1}(1-P)a_{\ell+1}'W^{q_\ell'}a_\ell'\cdots W^{q_1'}a_1')=0
\end{equation}
for all integers $k,\ell\ge0$ and $q_j,q_j'\ge1$ and all $a_ja_j'\in A$.
We will show~\eqref{eqn:psiaw} by induction on $k+\ell$.
For $k+\ell=0$, clearly $a_1(1-P)a_1'=a_1l^2(L^*)^2a_1'$ is in the kernel of $\Ec$,
hence also of $\psi$.
Suppose now $k+\ell>0$.
Writing $a_j=(a_j-\phi(a_j))+\phi(a_j)$ and similarly for $a_j'$ and then distributing,
we may assume that each $a_j$ and each $a_j'$ lies either in $\ker\phi$ or in $D$.
For every $d=(b_1,b_2)\in D$ we have $Ld=\alpha(d)L$,
where $\alpha$ is the automorphism of $D$ given by $\alpha\bigl((b_1,b_2)\bigr)=(b_2,b_1)$.
Indeed, we may take $d=d^*$ and then
\begin{eqnarray*}
(Ld-\alpha(d)L)^* (Ld - \alpha(d)L) & = & dL^* Ld -dL^* \alpha(d)L -L^* \alpha(d)
  Ld +L^* \alpha(d)\alpha(d)L \\
&=& d^2 -d^2 -d^2 +d^2 =0
\end{eqnarray*}
becase $L^* d L=\eta (d)=\alpha(d)$ for $d\in D$.
From $Ld=\alpha(d)L$ it follows that
\begin{equation*}
L^*d=\alpha(d)L^*,\quad Pd=dP\quad\text{and}\quad Wd=\alpha(d)W.
\end{equation*}
Note also that $W(1-P)=0$.
If $k\ge1$ and $a_{k+1}\in D$ then $W^{q_k}a_{k+1}(1-P)=W^{q_k}(1-P)a_{k+1}=0$;
similarly if $\ell\ge1$ and $a_{\ell+1}'\in D$ then $(1-P)a_{\ell+1}'W^{q_\ell'}=0$.
If $a_j\in D$ for some $1<j\le k$ then $W^{q_{j-1}}a_jW^{q_j}=W^{q_{j-1}+q_j}\alpha^{q_j}(a_j)$ and we may
use the induction hypothesis to conclude that~\eqref{eqn:psiaw} holds;
similarly if $a_j'\in D$ for some $1<j\le\ell$ then~\eqref{eqn:psiaw} holds.
Hence we may assume that $a_j\in\ker\phi$ for every $1<j\le k+1$ and
$a_j'\in\ker\phi$ for every $1<j\le\ell+1$.
Writing $W^{q_j}=\bigl(W^{q_j}-\psi(W^{q_j})\bigr)+\psi(W^{q_j})$ and similarly for $W^{q_j'}$, distributing
and letting $y_j=W^{q_j}-\psi(W^{q_j})$ and $y_j'=W^{q_j'}-\psi(W^{q_j'})$,
we find that $a_1W^{q_1}\cdots a_kW^{q_k}a_{k+1}(1-P)a_{\ell+1}'W^{q_\ell'}a_\ell'\cdots W^{q_1'}a_1'$
is equal to a sum of $2^{k+\ell}$ terms which are obtained by replacing each $W^{q_j}$
variously with $y_j$ and with $\psi(W^{q_j})$, and each $W^{q_j'}$
variously with $y_j'$ and with $\psi(W^{q_j'})$.
If  $z$ is one of these terms where where at least one $W^{q_j}$ or $W^{q_j'}$ has been replaced by its
expectation under $\psi$ then we can see that $\psi(z)=0$ by  using the induction hypothesis;
indeed, we write $y_j=W^{q_j}-\psi(W^{q_j})$ for each $y_j$ appearing in $z$ and similarly for each $y_j'$ and
then we distribute;
this expresses $z$ as a sum of terms to each of which the induction hypothesis applies to show has expectation zero under $\psi$.
We are left to show only that
\begin{equation*}
\Psi(a_1y_1\cdots a_ky_ka_{k+1}(1-P)a_{\ell+1}'y_\ell'a_\ell'\cdots y_1'a_1')=0.
\end{equation*}
But this holds by the freeness proved in Claim~\ref{claim:LAfree}.
Thus the proof of Claim~\ref{claim:Wunitary} is finished.
\end{proof}

\begin{claim}
The restriction of $\pi$ to $A$ is faithful.
\end{claim}
\begin{proof}
This follows from the fact that $\phi_1$ and $\phi_2$ have faithful GNS representations.
\end{proof}

The inner product arising from the GNS construction for $\psi$ gives a map $\psit:\pi(\Afr)\to D$ which, upon
identifying $D$ with $\pi(D)$ becomes a conditional expectation $\psit:\pi(\Afr)\to\pi(D)$;
moreover, we have $\psit\circ\pi=\pi\circ\psi{\restriction}_\Afr$.
\begin{claim}
\label{claim:freeM}
Let $q=\pi\bigl((1,0)\bigr)\in\pi(D)$ and let $V=\pi(W)$.
Then
\renewcommand{\labelenumi}{(\alph{enumi})}
\begin{enumerate}
\item $V$ is a unitary satisfying $V=V^*$ and $VqV=1-q$.
\end{enumerate}
Consider the subalgebras $M_1=q\pi(A)q$, $M_2=qV\pi(A)Vq=V(1-q)\pi(A)(1-q)V$
and $N=q\pi(D)q$.
Then
\renewcommand{\labelenumi}{(\alph{enumi})}
\begin{enumerate}
\addtocounter{enumi}{1}
\item $N\subseteq M_\iota$, ($\iota=1,2$);
\item $M_\iota$ is isomorphic to $A_\iota$ via an isomorphism that sends $N$ to $B$
      and conjugates $\psit{\restriction}_{M_\iota}$ to $\phi_\iota$, ($\iota=1,2$);
\item $M_1$ and $M_2$ are free with respect to $\psit$.
\end{enumerate}
\end{claim}
\begin{proof}
(a) follows from Claim~\ref{claim:PW}, Claim~\ref{claim:Wunitary}
and the fact that $W(b_1,b_2)=(b_2,b_1)W$ for all $(b_1,b_2)\in D$;
note that~(b) holds for the same reason.

For~(c), the isomorphisms are
\begin{align*}
M_1\ni\pi\bigl((a_1,0)\bigr)&\mapsto a_1\in A_1 \\
M_2\ni V\pi\bigl((0,a_2)\bigr)V&\mapsto a_2\in A_2,
\end{align*}
which we denote $\sigma_1$ and $\sigma_2$, respectively.
That $\sigma_1$ sends $N$ to $B$ and conjugates $\psit{\restriction}_{M_\iota}$ to $\phi_1$ is straightforward
to see.
We have $N\subseteq M_2$ because $\pi\bigl((b,0)\bigr)=V\pi\bigl((0,b)\bigr)V$ for every $b\in B$,
and this also shows that $N$ is mapped by $\sigma_2$ onto $B$.
We must show that
\begin{equation}
\label{eqn:sig2}
\sigma_2\circ\psit\bigl(V\pi\bigl((0,a_2)\bigr)V\bigr)=\phi_2(a_2)
\end{equation}
for every $a_2\in A_2$.
But
\begin{equation*}
\psit\bigl(V\pi\bigl((0,a_2)\bigr)V\bigr)=\psit\circ\pi\bigl(W(0,a_2)W\bigr)
=\pi\circ\psi\bigl(W(0,a_2)W\bigr).
\end{equation*}
Write $a_2=(a_2-\phi_2(a_2))+\phi_2(a_2)$.
It is easily seen using~\eqref{eqn:WL} that $\psi(W)=0$.
This and the freeness result proved in Claim~\ref{claim:LAfree} show that
$\psi\bigl(W(0,a_2-\phi_2(a_2))W\bigr)=0$,
while since $\phi_2(a_2)\in B$ we have
\begin{equation*}
\pi\circ\psi\bigl(W(0,\phi_2(a_2))W\bigr)=
\pi\circ\psi\bigl((\phi_2(a_2),0)\bigr)=
\pi\bigl((\phi_2(a_2),0)\bigr)=
V\pi\bigl((0,\phi_2(a_2))\bigr)V
\overset{\sigma_2}\mapsto\phi_2(a_2).
\end{equation*}
Hence~\eqref{eqn:sig2} holds and~(c) is proved.

To prove~(d) it will suffice to show that $\psit(x)=0$ whenever
$x=x_1x_2\cdots x_n$ where $n\ge1$, $x_j\in M_{\iota_j}\cap\ker\psit$
and $\iota_1\ne\iota_2,\,\iota_2\ne\iota_3,\ldots,\iota_{n-1}\ne\iota_n$.
Writing
\begin{equation*}
x_j=\begin{cases}
\pi\bigl((a_j,0)\bigr)\text{ for }a_j\in A_1\cap\ker\phi_1&\text{if }\iota_j=1 \\
\pi\bigl(W(0,a_j)W\bigr)\text{ for }a_j\in A_2\cap\ker\phi_2&\text{if }\iota_j=2,
\end{cases}
\end{equation*}
we find that $\psit(x)=\psit\circ\pi(y)=\pi\circ\psit(y)$ where $y=y_1y_2\cdots y_n$ and
\begin{equation*}
y_j=\begin{cases}
(a_j,0)&\text{if }\iota_j=1 \\
(0,a_j)&\text{if }\iota_j=2.
\end{cases}
\end{equation*}
Rewriting $y$ as a product of $W$'s alternating with $(0,a_j)$'s and $(a_j,0)$'s, using that
each $(a_j,0)$ and $(0,a_j)$ lies in $\ker\psi$, that $\psi(W)=0$ and the freeness proved in
Claim~\ref{claim:LAfree}, we find that $\psi(y)=0$.
This finishes the proof of Claim~\ref{claim:freeM}.
\end{proof}

We now finish the proof of Proposition~\ref{prop:AtemEH}.
By properties of the free product construction,
$\At$ is isomorphic to the image of $C^*(M_1\cup M_2)$ under the GNS representation of
$\psit{\restriction}_{C^*(M_1\cup M_2)}$,
while $C^*(M_1\cup M_2)$ is itself a subalgebra of a quotient of a subalgebra of $E(H)$.
\end{proof}

\begin{cor}[\cite{D:ExactFP}]
\label{cor:AmFPEx}
Let $B$ be a unital C$^*$--algebra and let $A_1$ and $A_2$ be unital C$^*$--algebras each containing
a copy of $B$ as a unital C$^*$--subalgebra and having a conditional expectation
$\phi_\iota:A_\iota\to B$ whose GNS representation is faithful.
Let
\begin{equation*}
(\At,\phit)=(A_1,\phi_1)*_B(A_2,\phi_2)
\end{equation*}
be the reduced amalgamated free product of C$^*$--algebras.
Then $\At$ is an exact C$^*$--algebra if and only if $A_1$ and $A_2$ are exact.
\end{cor}
\begin{proof}
Since $A_1$ and $A_2$ are C$^*$--subalgebras of $\At$, exactness of $\At$ implies
that of $A_1$ and $A_2$.

For the converse, suppose $A_1$ and $A_2$ are exact and let $A=A_1\oplus A_2$.
Then $A$ is an exact C$^*$--algebra.
By Theorem~\ref{thm:EHexact}, $E(H)$ is an exact C$^*$--algebra whenever $H$ is a Hilbert $A,A$--bimodule.
Using Proposition~\ref{prop:AtemEH}, the well known fact that C$^*$--subalgebras of exact C$^*$--algebras are exact
and Kirchberg's result~\cite{Kirchberg:ComUHF} that exactness passes to quotients, we have that $\At$ is exact.
\end{proof}

\section{Bogljubov automorphisms.}
\label{sec:Bog}

In this section we consider the analogues of Bogljubov automorphisms on the extended
Cuntz--Pimnser algebras $E(H)$ when $H$ is a $B,B$--bimodule and we prove that if $B$ is
finite dimensional then the topological entropy of every Bogljubov automorphism is zero.
We begin, however, by showing that if $B$ is a finite dimensional C$^*$--algebra
and if $H$ is a Hilbert $B$--module then every countably generated Hilbert $B$--submodule
of $H$ is complemented in $H$.
We are convinced this is well known, but we don't know of a reference.

\begin{prop}[Gram--Schmidt procedure]
\label{prop:G-S}
Let $B$ be a finite dimensional C$^*$--algebra and let $H$ be a right Hilbert $B$--module.
Let $X$ be a finite or countably infinite subset of $H$.
Then there is a finite or countably infinite subset $V$ of $H$ such that
\renewcommand{\labelenumi}{(\roman{enumi})}
\begin{enumerate}
\item the submodule of $H$ generated by $V$ equals the submodule of $H$ generated by $X$,
\item if $v,w\in V$ and $v\ne w$ then $\langle v,w\rangle=0$,
\item if $v\in V$ then $\langle v,v\rangle$ is a minimal projection in $B$.
\end{enumerate}
\end{prop}

\begin{proof}
Let $e_1,\ldots,e_n$ be minimal projections in $B$ such that $e_1+e_2+\cdots+e_n=1$.
We may without loss of generality assume that for every $x\in X$, $x=xe_\ell$
for some $\ell\in\{1,\ldots,n\}$.
Let $X=\{x_1,x_2,\ldots\}$ be an enumeration of $X$ and let $\ell_j$ be such that
$x_j=x_je_{\ell_j}$.
Note that $\langle x_1,x_1\rangle=\nm{x_1}^2e_{\ell_1}$ and let
$$ v_1=
\begin{cases}
0&\text{if }x_1=0 \\
\nm{x_1}^{-1}x_1&\text{if }x_1\ne0.
\end{cases} $$
Then either $v_1=0$ or $\langle v_1,v_1\rangle=e_{\ell_1}$.
We now recursively define elements $v_2,v_3,\ldots$ so that for all $j$ we have
$v_je_{\ell_j}=v_j$, $\langle v_j,v_j\rangle\in\{0,e_{\ell_j}\}$
and if $i<j$ then $\langle v_i,v_j\rangle=0$.
For the recursive step, if $n\ge2$, if $v_1,\ldots, v_{n-1}$ have been defined
and if $X$ has at least $n$ elements then let
$w_n=x_n-\sum_{j=1}^{n-1}v_j\langle v_j,x_n\rangle$.
Then for every $i\in\{1,\ldots,n-1\}$ we have
$\langle v_i,w_n\rangle=\langle v_i,x_n\rangle-\bigl\langle v_i,v_i\langle v_i,x_n\rangle\bigr\rangle
=\langle v_i,x_n\rangle-\langle v_i,v_i\rangle\langle v_i,x_n\rangle=0$.
Note that $w_n=w_ne_{\ell_n}$;
let
$$ v_n=
\begin{cases}
0&\text{if }w_n=0 \\
\nm{w_n}^{-1}w_n&\text{if }w_n\ne0.
\end{cases} $$
Letting $V=\{v_j\mid v_j\ne0\}$ does the job.
\end{proof}

Although we will only apply the following proposition when the submodule $K$ is assumed
to be finitely generated, for sake of completeness we would like to give the more general result.

\begin{prop}
\label{prop:compl}
Let $B$ be a finite dimensional C$^*$--algebra, let $H$ be a right Hilbert $B$--module
and let $K$ be a closed $B$--submodule of $H$ that is finitely or countably generated,
(meaning that $K$ has a dense submodule that is finitely or countably generated).
Then $K$ is a complemented submodule of $H$.
\end{prop}

\begin{proof}
Let $X$ be a finite or countable set such that the submodule of $H$ generated by $X$
is dense in $K$.
Let $V$ be the set obtained form $X$ using the Gram--Schmidt procedure of Propositon~\ref{prop:G-S}.
We shall define $P:H\to H$ by
$$ Ph=\sum_{v\in V}v\langle v,h\rangle, $$
where when $V$ is infinite we shall show that the sum converges in $H$.
Suppose first that $V$ is finite.
Then easy calculations show that $P\in\LEu(H)$, $P^2=P=P^*$ and consequently $\nm P\le1$.

Now suppose that $V$ is infinite and enumerate it by $V=\{v_1,v_2,\ldots\}$.
For every positive integer $n$ and $h\in H$ let $P_nh=\sum_{j=1}^nv_j\langle v_j,h\rangle$.
Then by the result for the case of $V$ finite we have
\begin{equation*}
\langle h,h\rangle\ge\langle P_nh,P_nh\rangle
=\sum_{j=1}^n\langle h,v_j\rangle\langle v_j,h\rangle.
\end{equation*}
For every state $\phi$ on $B$, the sequence
$\Bigl(\phi\bigl(\langle P_nh,P_nh\rangle\bigr)\Bigr)_{n=1}^\infty$ is a bounded and increasing
sequence of positive numbers, hence converges.
Therefore the sequence $\bigl(\langle P_nh,P_nh\rangle\bigr)_{n=1}^\infty$ converges in $B$.
Then for $m<n$ we have that
\begin{equation*}
\langle P_nh-P_mh,P_nh-P_mh\rangle=\sum_{j=m+1}^n\langle h,v_j\rangle\langle v_j,h\rangle
\le\sum_{j=m+1}^\infty\langle h,v_j\rangle\langle v_j,h\rangle
\end{equation*}
and the right--hand--side tends to zero as $m\to\infty$.
Therefore the sequence $(P_nh)_{n=1}^\infty$ is Cauchy in $H$, hence converges in $H$.
We may thus define $Ph=\sum_{j=1}^\infty v_j\langle v_j,h\rangle$.
From the corresponding facts for the finite dimensional case
we obtain that $P\in\LEu(H)$, $P^2=P=P^*$ and $\nm P\le1$.

It remains to show that $PH=K$.
If $V$ is finite then this is clear, so suppose $V$ is infinite.
Given $h\in H$ we have $Ph=\lim_{n\to\infty}P_nh$ and $P_nh\in K$, so $Ph\in K$.
Given $k\in K$ then for all $\epsilon>0$ there are $n\ge1$ and
$k_\epsilon\in\lspan\{v_1,\ldots,v_n\}$ such that $\nm{k-k_\epsilon}<\epsilon$.
But $Pk_\epsilon=k_\epsilon$ and hence $\nm{Pk-k}\le2\epsilon+\nm{Pk_\epsilon-k_\epsilon}=2\epsilon$.
Therefore $Pk=k$.
\end{proof}

\begin{defi}
\label{def:Bog}
Let $B$ be a C$^*$--algebra and let $H$ be a Hilbert $B,B$--bimodule such that
$\{\langle h_1,h_2\rangle\mid h_1,h_2\in H\}$ generates $B$;
suppose that $U:H\to H$ is a $\Cpx$--linear map such that
for some automorphism $\beta$ of $B$ we have
\begin{align*}
\bigl\langle U(h_1),U(h_2)\bigr\rangle=\beta(\langle h_1,h_2\rangle),&\quad h_1,h_2\in H \\
U(b_1hb_2)=\beta(b_1)U(h)\beta(b_2),&\quad h\in H,\,b_1,b_2\in B.
\end{align*}
(Note that $\beta$ is uniquely determined by $U$ and the first of the above equations.)
Then there is an automorphism $E(U)$ of $E(H)$, given by $E(U)(l(h))=l(Uh)$ ($h\in H$)
and $E(U)(b)=\beta(b)$ ($b\in B$).
We call $E(U)$ the \emph{Bogljubov automorphism} of $E(H)$ associated to $U$.
\end{defi}

\begin{thm}
\label{thm:Bog}
Let $B$ be a finite dimensional C$^*$--algebra and let $H$ be a Hilbert $B,B$--bimodule
such that $\{\langle h_1,h_2\rangle\mid h_1,h_2\in H\}$ generates $B$.
Then the topological entropy of every Bogljubov automorphism $E(U)$ of $E(H)$
is zero.
\end{thm}

\begin{proof}
The Bogljubov automorphism $E(U)$ arises from a $\Cpx$--linear map $U:H\to H$ and
an affiliated automorphism $\beta\in\Aut(B)$ as in Definition~\ref{def:Bog}.
Let $\Ht=H\oplus B$ and $\xi=0\oplus1\in\Ht$.
Let $\Ut:\Ht\to\Ht$ be defined by $\Ut(h\oplus b)=U(h)\oplus\beta(b)$, ($h\in H,\,b\in B$);
note that $\Ut$ and $\beta$ together satisfy the conditions in Definition~\ref{def:Bog}, so that
we have the Bogljubov automorphism $E(\Ut)$ of $E(\Ht)$.
Now $E(H)$ is canonically embedded in $E(\Ht)$ and the restriction of $E(\Ut)$ to $E(H)$ is $E(U)$;
by the monotonicity of $ht$, which was proved in~\cite[Proposition 2.1]{Brown:extopent},
it will therefore suffice to show that $ht(E(\Ut))=0$.

As in the proof of Theorem~\ref{thm:EHexact}, $E(\Ht)$ is isomorphic to
the crossed product C$^*$--algebra $A\rtimes_\Psi\Nats$, where $A=\clspan\Omega$ with
$$ \Omega=B\cup\{l(h_1)\cdots l(h_m)l(h_{m+1})^*\cdots l(h_{2m})^*\mid m\ge1,\,h_j\in\Ht\} $$
and where $\Psi$ is the endomorphism of $A$ given by $\Psi(x)=LxL^*$ with $L=l(\xi)$.
Since $\Ut(\xi)=\xi$, we have $E(\Ut)L=L$ and hence the restiction of $E(\Ut)$ to $A$
is an automorphism of $A$ that commutes with the endomorphism $\Psi$.
It thus follows from Proposition~\ref{prop:htbetahatend}
that $ht(E(\Ut))=ht(E(\Ut){\restriction}_A)$.
Let $\gamma$ denote the automorphism $E(\Ut){\restriction}_A$ of $A$.

Let $\tau:B\to\LEu(\VEu)$ be a faithful unital representation of $B$ on a Hilbert space $\VEu$
and let $\pi$ denote the representation of $A$ which is
the inclusion $A\hookrightarrow\LEu(\Fc(\Ht))$ followed by
the representation $x\mapsto x\otimes1$ of $\LEu(\Fc(\Ht))$ on the Hilbert space
$\Fc(\Ht)\otimes_\tau\VEu$.
We must show $ht(\gamma)=0$ and to do so it will suffice to show
that $ht(\pi,\gamma,\omega,\delta)=0$ for every $\delta>0$
and every finite subset $\omega$ of $\Omega$.
Given a finite subset $\omega\subseteq\Omega$ there are $n\in\Nats$
and a $B,B$--subbimodule $K$ of $\Ht$ that is finite dimensional as a $\Cpx$--vector space
and such that $\omega\subseteq\Omega(n,K)$ where
$$ \Omega(n,K)\eqdef B\cup
\{l(h_1)\cdots l(h_m)l(h_{m+1})^*\cdots l(h_{2m})^*\mid1\le m\le n,\,h_j\in K\}. $$
Recall the definition of $\Fc_n(\Ht)$ from the proof of Claim~\ref{claim:Aexact} and
let $P_n\in\LEu(\Fc(\Ht))$ denote the projection onto $\Fc_n(\Ht)$.
Consider the completely positive contractions
\begin{alignat*}{2}
\Phi_n:&\LEu(\Fc(\Ht))\to\LEu(\Fc_n(\Ht)),&\qquad\Phi_n(x)&=P_nxP_n \\
\Psi_n:&\LEu(\Fc_n(\Ht))\to\LEu(\Fc(\Ht)),&\qquad\Psi_n(y)&=W_n^*(y\otimes1)W_n,
\end{alignat*}
where $W_n:\Fc(\Ht)\to\Fc_n(\Ht)\otimes_B\Fc(\Ht^{\otimes_B(n+1)})$ is the unitary operator
canonically defined by the decomposition in
equation~(\ref{eqn:Fcdecomp}) in the proof of Claim~\ref{claim:Aexact}.
Note that $\Psi_n\circ\Phi_n(x)=x$ for every $x\in\Omega(n,\Ht)$.

For every integer $p\ge1$ let $K_p=K+\Ut(K)+\Ut^2(K)+\cdots+\Ut^{p-1}(K)$;
then
$$ \omega\cup\gamma(\omega)\cup\cdots\cup\gamma^{p-1}(\omega)\subseteq\Omega(n,K_p). $$
Moreover, $K_p$ is an $B,B$--subbimodule of $\Ht$ whose $\Cpx$--linear dimension satisfies
$\dim_\Cpx(K_p)\le p\dim_\Cpx(K)$.
Let
$$ \Fc_n(K_p)=B\oplus\bigoplus_{k=1}^nK_p^{(\otimes_B)k}\subseteq\Fc_n(\Ht)\subseteq\Fc(\Ht). $$
Clearly $\Fc_n(K_p)$ is a finite dimensional $B,B$--subbimodule of $\Fc_n(\Ht)$.
By Proposition~\ref{prop:compl}, there is $Q_{n,p}\in\LEu(\Fc_n(\Ht))$ which is the projection
onto $\Fc_n(K_p)$;
note that $Q_{n,p}$ commutes with the left action of $B$ on $\Fc_n(\Ht)$.
Consider the completely positive contractions
\begin{alignat*}{2}
\Theta_{n,p}:&\LEu(\Fc_n(\Ht))\to\LEu(\Fc_n(K_p)),&\qquad\Theta_{n,p}(x)&=Q_{n,p}xQ_{n,p} \\
\Upsilon_{n,p}:&\LEu(\Fc_n(K_p))\to\LEu(\Fc_n(\Ht)),&
 \qquad\Upsilon_{n,p}(y)&=Q_{n,p}yQ_{n,p}+V_{B,n}^*yV_{B,n}(1-Q_{n,p})
\end{alignat*}
where $V_{B,n}\in\LEu(B,\Fc_n(\Ht))$ maps $B$ to the submodule $B\oplus0$ of $\Fc_n(\Ht)$
via $b\mapsto b\oplus0$.
Since $l(h)^*(1-Q_{n,p})=0$ whenever $h\in K_p$, we see that
$\Psi_n\circ\Upsilon_{n,p}\circ\Theta_{n,p}\circ\Phi_n(x)=x$
for every $x\in\Omega(n,K_p)$.
As $\LEu(\Fc_n(K_p))$ is a finite dimensional C$^*$--algebra, we have that
$$ rcp(\pi,\omega\cup\gamma(\omega)\cup\cdots\cup\gamma^{p-1}(\omega),\delta)
\le\rank\LEu(\Fc_n(K_p)). $$
As the C$^*$--algebra $\LEu(\Fc_n(K_p))$ can be faithfully represented on the Hilbert space
$\Fc_n(K_p)\otimes_B\VEu$, making a crude estimate we get
\begin{align*}
\rank\LEu(\Fc_n(K_p))&\le\dim\bigl(\Fc_n(K_p)\otimes_B\VEu\bigr)
\le\dim(\VEu)\Bigl(\sum_{k=0}^n\dim_\Cpx(K_p)^k\Bigr) \\
&\le n\dim(\VEu)\dim_\Cpx(K_p)^n\le n p^n\dim(\VEu)\dim_\Cpx(K)^n.
\end{align*}
Because this upper bound grows subexponentially as $p\to\infty$, we conclude that
$ht(\pi,\gamma,\omega,\delta)=0$.
\end{proof}

\section*{Acknowledgements}
The authors would like to thank the Erwin Schr\"odinger Institute (Vienna) and J.B.\ Cooper
for organizing the ESI's program in Functional Analysis in spring 1999,
where our collaboration on the subject of this paper began.
Part of this research was done while the first named author was empolyed at Odense University
and was partially supported by a grant from the Danish National Research Foundation.
This work was completed while the second named author was
participanting in the program  on Free Probability and Operator Spaces at the
Institute Henri Poincar\'e (Paris).
Finally, we would like to thank \'E.\ Blanchard for pointing out a mistake
in equation~\eqref{eqn:Fcdecomp} of our original version.

\end{spacing}

\newpage

\begin{spacing}{1.0}

\bibliographystyle{plain}

\vskip1ex
\scriptsize
\noindent{\sc Department of Mathamtics, Texas A\&M University, College Station TX 77843--3368, USA} \\
\noindent{\sl E-mail:} {\tt Ken.Dykema@math.tamu.edu} \\
\noindent{\sl Internet URL:} {\tt http://www.math.tamu.edu/$\,\widetilde{\;}$Ken.Dykema/} \\
\vskip1ex
\noindent{\sc Department of Mathamtics, University of California, 405 Hilgard Ave., Los Angeles CA 90095--1555, USA} \\
\noindent{\sl E-mail:} {\tt shlyakht@math.ucla.edu} \\

\end{spacing}
\end{document}